\title{ ~~\\ On some claims in Ramanujan's `unpublished' manuscript on the 
partition and tau functions}
\author{Pieter Moree}
\documentclass[12pt]{article}
\usepackage{amssymb, latexsym, amsfonts}
\def\@ptsize{2}
\newtheorem{Thm}{Theorem}

\newtheorem{Lem}{Lemma}
\newtheorem{Cor}{Corollary}

\newcommand{\qed}{\hfill $\Box$}

\begin{document}
\date{}
\maketitle
{\def\thefootnote{}
\footnote{\noindent P. Moree: KdV Institute, University of Amsterdam,
Plantage Muidergracht 24, 1018 TV Amsterdam, The Netherlands, e-mail:
moree@science.uva.nl}}
{\def\thefootnote{}
\footnote{{\it Mathematics Subject Classification (2001)}.
01A99, 11N37, 11Y60}}
\begin{abstract}
\noindent We disprove various claims made by Ramanujan in his, until
very recently, unpublished manuscript \cite{bono} on the partition and
tau functions.
Furthermore, the second part of a related paper by G.K. Stanley \cite{stanley} 
is corrected (the first part of which was earlier
corrected by D. Shanks \cite{shanks}).
\end{abstract}
\section{Introduction}
In his first letter (16 Jan. 1913) to Hardy \cite[p. 24]{bera}, Ramanujan made various claims. The
fourth of them
reads as follows:
"$(4)~~1,2,4,5,8,9,10,13,16,17,18,\cdots$ are numbers which are
either themselves squares or which can be
expressed as the sum of two squares. The number of such numbers greater than $A$ and 
less than $B$ 
\begin{equation}
\label{nulde}
=K\int_A^B{dx\over \sqrt{\log x}}+\theta(x),
\end{equation}
where $K=0.764\cdots$ and $\theta(x)$ is very small when compared with the previous integral. $K$ and
$\theta(x)$ have been exactly found, though complicated....". (Note that $\theta(x)$ should
be $\theta(B)$.)
Answering an inquiry of Hardy \cite[p. 49]{bera}, Ramanujan in his second letter to Hardy 
\cite[p. 56]{bera} claims:
"The order of $\theta(x)$ which you asked in your letter is $\sqrt{x}/\sqrt{\log x}$". In his
lectures \cite{hardy} on Ramanujan's work, Hardy states that Ramanujan also gave the exact value of $K$, namely
$$K={1\over \sqrt{2}}\prod_{p\equiv 3({\rm mod~}4)}\left(1-{1\over p^2}\right)^{-1/2}.$$
Note that Ramanujan's claim, if correct, would imply that
$$B(x)={Kx\over \sqrt{\log x}}\left(1+{C_2\over \log x}+O\left({1\over \log ^2 x}\right)\right),$$
with $C_2=1/2$, where $B(x)$ denotes the number of integers $n\le x$
of the form $u^2+v^2$ with $u$ and $v$ integers.
Landau proved in 1908, using contour integration, that asymptotically
$B(x)\sim Kx/\sqrt{\log x}$. His method can be easily
extended \cite{serre} to prove that $B(x)$ has an asymptotic expansion in the sense
of Poincar\'e, namely for every integer $r\ge 2$, we have
$$B(x)={Kx\over \sqrt{\log x}}\left(1+\sum_{j=2}^{r}
{C_j\over \log ^{j-1}x}+O\left({1\over \log ^r x}\right)\right),$$
with $C_2,\cdots,C_r$ constants.
Shanks \cite{shanks}, correcting on Stanley \cite{stanley}, computed the
second order constant $C_2$ to equal $0.5819486\cdots$ and thus disproved Ramanujan's claim. For an overview of further 
results in this direction the reader is referred to \cite{yana}.\\
\indent There is some evidence 
(but see \cite[p. 92]{rankinzelf}) that along with his final letter (12 Jan. 1920) to Hardy,
Ramanujan included a manuscript on congruence properties of $\tau(n)$ and $p(n)$. 
In this manuscript Ramanujan considers, for various special small primes
$q$, the quantity $\sum_{n\le x,~q\nmid \tau(n)}1$ and makes claims
similar to (\ref{nulde}). 
He defines $t_n=1$ if $\tau(n)\not\equiv 0 ({\rm mod~}q)$ and
$t_n=0$ otherwise. He then typically writes: "It is easy to prove by quite
elementary methods that $\sum_{k=1}^n t_k=o(n)$. It can be shown by
transcendental methods that
\begin{equation}
\label{ditgaatnoggoed}
\sum_{k=1}^n t_k\sim {Cn\over (\log n)^{\delta}};
\end{equation}
and
\begin{equation}
\label{valseanalogie}
\sum_{k=1}^n t_k=C\int_1^n {dx\over (\log x)^{\delta}}+O\left({n\over (\log n)^r}
\right).
\end{equation}
where $r$ is any positive number". 
Note that the truth of $\sum_{k=1}^n t_k=o(n)$
would imply that $q|\tau(n)$ for almost all $n$.
The values of $\delta$ can be found in the final
column of Table 1.
Except for $q=5$ and $q=691$ Ramanujan also writes down an Euler product
for $C$. It is not difficult to check that these are correct, except when $q=23$, in
which case this is due to a factor $(1-23^{-s})^{-1}$ erroneously omitted in the
generating function.\\ 
\indent It appears from \cite{stanley} that Hardy planned to have this manuscript published 
under Ramanujan's name after some editing. Indeed, he 
published some parts of it (see \cite{bono}). Unfortunately, he never published Ramanujan's
full manuscript (which indeed needed some serious editing). Some of his results
in connection with the unpublished
manuscript were further worked out by his research student Geraldine Stanley and
published in 1928 \cite{stanley}. She claimed Ramanujan's assertion regarding
$C_2$ and the second order coefficient for $\sum_{n\le x,~5\nmid \tau(n)}1$ to be false.
Unfortunately, her paper contains several typos and
some mathematical errors, which were corrected by Shanks \cite{shanks} for the $B(x)$ case
and are corrected in this paper for the $5|\tau(n)$ case in Section \ref{Geraldine}.\\
\indent In 1928 Hardy passed on his Ramanujan materials to G.N. Watson. 
Watson's papers \cite{watson1, watson2, watson3} are inspired by Ramanujan's `unpublished' manuscript.
In particular Watson \cite{watson1} rigorously proved 
Ramanujan's assertion that $691|\tau(n)$ for almost all $n$
by establishing (\ref{ditgaatnoggoed}) with $q=691$ and $\delta=1/690$. 
A remarkable result, since
for $1\le n\le 5000$ we have, as Ramanujan computed, $691\nmid \tau(n)$, 
except for
the multiples of $1381$ in that range.
After Watson's death in 1965 a part of the manuscript came into the possession of
the library of Trinity College, Cambridge. Watson's copy of the remaining 
part can be found in the library of Oxford's Mathematical Institute.
A full version of the unpublished manuscript (with proofs
and commentary) has only recently become available to the
general public \cite{bono}.\\
\indent Ramanujan's claims are only the tip of an iceberg: it is now known 
\cite{serre, serre2} that if
$\sum_{n=1}^{\infty}a_nq^n$, with $q=e^{2\pi i z}$, is the Fourier expansion of a modular
form of integral weight with integral coefficients for a congruence
subgroup of the full modular group, then for every positive integer $M$,
for almost all $n$ we have $a(n)\equiv 0({\rm mod~}M)$. In particular, we have
that (\ref{ditgaatnoggoed}) holds true for every odd prime $q$ not in Table 1 with 
$\delta=q/(q^2-1)$. It can also be shown \cite{SirSWD}, using
$l$-adic representations, that $2$ and the primes
in Table 1 are the only primes for which congruences for $\tau(n)$ exist.\\
\indent The purpose of this paper is to establish the 
falsity of the various claims made
by Ramanujan of the
type (\ref{valseanalogie}). 
We do this by evaluating the relevant second order coefficients and
compute them with at least a few decimals of precision.
Moreover, we will correct the part of Stanley's paper pertaining to the
$5|\tau(n)$ case. As far as possible we adopt Ramanujan's notation and the order
in which we deal with the various primes follows the order in which
they appear in Ramanujan's manuscript.

\section{Some relevant general results}
In this section we quote some general results that allow us to evaluate the second-order
constants of the functions considered by Ramanujan.\\
\indent The following result was obtained by Moree \cite{moree} using elementary methods. It
can be deduced on analysing an appropriate functional equation 
(apparently first considered by Levin and Fainleib) for the relevant counting function.
Given a multiplicative function $f$, we define $\Lambda_f(n)$ implicitly by
$$f(n)\log n=\sum_{d|n}f(d)\Lambda_f({n\over d}).$$
We denote the formal Dirichlet series $\sum_{n=1}^{\infty}f(n)n^{-s}$ associated
to $f$ by $L_f(s)$. It is easy to see that
 $-L_f'(s)/L_f(s)=\sum_{n=1}^{\infty}\Lambda_f(n)n^{-s}$. As usual the logarithmic integral
 $\int_2^x{dt/\log t}$ is denoted by Li$(x)$.
\begin{Thm} {\rm \cite{moree}}.
\label{oud}
Let $f$ be a multiplicative function
satisfying
\begin{equation}
\label{grens}
0\le f(p^r)\le c_1c_2^r,~c_1\ge 1,~1\le c_2<2,
\end{equation}
and
$\sum_{p\le x}f(p)=\tau {\rm Li}(x)+O\left({x\log^{-2-\rho}x}\right),$
where $\tau$ and $\rho>1$ are positive real fixed numbers. Then, for
some $\epsilon>0$ and constant $B_f$ we have
\begin{equation}
\label{nognodig}
\sum_{n\le x}{\Lambda_f(n)\over n}=\tau \log x + B_f
+O\left({x\over \log^{1+\epsilon} x}\right).
\end{equation}
There exist constants $\epsilon'>0$ and constants $C_1(f)>0$ and $C_2(f)$ such that
$$\sum_{n\le x}f(n)={C_1(f)x\over \log ^{1-\tau} x}
\left(1+{C_2(f)\over \log x}+O\left({1\over \log ^{1+\epsilon'} x}\right)\right),
~x\rightarrow\infty.$$
In particular, $C_2(f)=(1-\tau)(1+B_f)$.
\end{Thm}
Remark. Alternatively we can write
$$C_2(f)=(1-\tau)\left(1+\tau \gamma+\sum_{n=1}^{\infty}{\Lambda_f(n)-\tau\over
n}\right),$$
where $\gamma$ denotes Euler's constant.\\
\begin{Cor}
\label{bisnul}
Let $f$ satisfy the conditions of Theorem {\rm \ref{oud}} with $\tau\ne 1$. Suppose that
$$
\sum_{n\le x}f(n)=C_1(f)\int_1^x{dt\over \log^{1-\tau}t}+O\left({x\over \log^r x}\right),$$
for some $r>2-\tau$, then we must have $B_f=0$.
\end{Cor}
\indent The constant $B_f$ appearing in
(\ref{nognodig}) can be computed using the following
result.
\begin{Lem}
\label{edmundo}
Let $f$ satisfy the conditions of Theorem {\rm \ref{oud}} and let $B_f$ and
$\tau$ be defined as in Theorem {\rm \ref{oud}}. 
Then
$$B_f=
-\lim_{s\rightarrow 1+0}\left({L_f'(s)\over L_f(s)}+{\tau\over s-1}\right).$$
\end{Lem}
Example. Take $f={\bf 1}$, that 
is $f(n)=1$ for every $n\ge 1$. Then $\Lambda_f$ equals the Von Mangoldt function
and $L_f(s)=\zeta(s)$. Using the well-known Taylor expansion 
\begin{equation}
\label{stieltjes}
\zeta(s)={1\over s-1}+\gamma+o(1) 
\end{equation}
around
$s=1$, we obtain that $B_f=-\gamma$.\\
\indent Our proof of Lemma \ref{edmundo} will make use of the following result due to
Landau \cite[pp. 73-74]{landau}.
\begin{Lem}
\label{hulp}
Suppose that $\sum_{n\le x}h(n)=\alpha x+O(g(x)),$ where $g(x)$ is a positive
function of $x$ such that $g(x)/x^2$ is monotonically decreasing for every
$x\ge x_0$, $x_0$ some fixed number, and where
$\int_{x_0}^{\infty}g(t)dt/t^2<\infty$. Then, for Re$(s)>1$, we have the
estimate
$$\sum_{n=1}^{\infty}{h(n)\over n^s}={\alpha\over s-1}+
\beta+o(1),~s\rightarrow 1+0,$$
for some constant $\beta$. Furthermore, we have
$$\sum_{n\le x}{h(n)\over n}=\alpha \log x +\beta+o(1),~x\rightarrow\infty.$$
\end{Lem}
\indent We are now in the position to prove Lemma \ref{edmundo}.\\

\noindent {\it Proof of Lemma} \ref{edmundo}. By partial integration
we deduce from (\ref{nognodig}) that
$\sum_{n\le x}\Lambda_f(n)=\tau x + O(x\log^{-1-\epsilon}x)$. We thus can
invoke Lemma \ref{hulp} with $h(n)=\Lambda_f(n)$, $\alpha=\tau$ and
$g(x)=x\log^{-1-\epsilon}x$. On noting that
$-L_f'(s)/L_f(s)=\sum_{n=1}^{\infty}\Lambda_f(n)n^{-s}$, the result
then follows. \qed\\

\indent Serre \cite{serre} gave some beautiful applications of Landau's method to
counting functions involving coefficients of modular forms. In order to formulate
his result we have to define the concept of {\it Frobenius} density.
A set of primes $\cal P$ is called Frobenius of density $\delta$, if there exists
a finite Galois extension $K/\mathbb Q$ and a subset $H$ of $G:={\rm Gal}(K/\mathbb Q)$
such that $H$ is stable under conjugation, $|H|/|G|=\delta$ and for every prime $p$, with
at most finitely many exceptions, one has $p\in {\cal P}$ if and only if
$\sigma_p(K/\mathbb Q)\in H$, where $\sigma_p(K/\mathbb Q)$ denotes the Frobenius
map of $p$ in $G$ (defined modulo conjugation in case $p$ does not divide the
discriminant of $K$).
Serre \cite{serre}, using's Landau's contour integration method, established 
the following result, which we formulate here in a slightly less general
form.
\begin{Thm} {\rm \cite[Th\'eor\`eme 2.8]{serre}}.
\label{Serre}
Let $\alpha(n):\mathbb N\rightarrow \mathbb Z$ be a multiplicative function. Let $q$ be
a fixed prime. Suppose that the set $P:=\{p{\rm ~is~prime}:~q|\alpha(p)\}$ is Frobenian of
density $0<1-\tau<1$. Let $h(s)=\sum_{q\nmid \alpha(n)}n^{-s}$. Then around $s=1$ we have
$${h(s)\over s}={1\over (s-1)^{\tau}}(e_0+e_1(s-1)+\cdots+e_k(s-1)^k+\cdots),$$
for certain numbers $e_j$ with $e_0\ne 0$ and we have, furthermore,
$$\sum_{n\le x\atop q\nmid \alpha(n)}1={e_0x\over \Gamma(\tau)\log^{1-\tau}x}\left(1
+\cdots+{\Gamma(\tau)e_k\over \Gamma(\tau-k)e_0\log^k x}+O\left({1\over \log ^{k+1} x}\right)\right).$$
In particular the second order constant, $\Gamma(\tau)e_1/(\Gamma(\tau-1)e_0)$, equals
$$(1-\tau)\lim_{s\rightarrow 1+0}\left(-{L_f'(s)\over L_f(s)}-{\tau\over s-1}+1\right).$$
\end{Thm}
It can be shown that if the conditions of Theorem \ref{Serre}
are satisfied, also the conditions of Theorem \ref{oud} are satisfied.
Theorem \ref{oud} in combination with Lemma \ref{edmundo} and Theorem
\ref{Serre} both predict the same second order coefficient.\hfil\break
\indent Let
$T(s)=\sum_{n=1}^{\infty}t_n/n^s$.
The approach followed in disproving Ramanujan's claims of the format (\ref{valseanalogie})
is to write $T(s)=\zeta(s)^{\tau}g(s)$ with $g(s)$ a regular function for Re$(s)>1/2$. By
Lemma \ref{edmundo} and (\ref{stieltjes})  it then follows that $B_t=-\tau\gamma-g'(1)/g(1)$. 
We have $\tau=1-\delta$.
The numerical work (carried out in Section \ref{numerical}) shows that $B_t\ne 0$.
From Corollary \ref{bisnul} the falsity of Ramanujan's claim then follows
for every $r>2-\tau$.

\section{Divisibility of tau by 2}
Ramanujan shows that $\tau(n)$ is odd or even according as $n$ is an odd square or not.
It thus follows that $\sum_{n\le x,~2\nmid \tau(n)}1=[{1+\sqrt{x}\over 2}]$.

\section{Divisibility of tau by 5}
\label{five}
At \cite[p. 47]{bono} Ramanujan makes a statement of the form 
(\ref{valseanalogie}) with $q=5$ and $\delta=1/4$.
Put $t_n=0$ if $5|\tau(n)$ and $t_n=1$ otherwise. Let
$T(s)=\sum_{n=1}^{\infty}t_n/n^s$. 
Denote $\sum_{d|n}d^r$ by $\sigma_r(n)$.
On using that $\tau(n)\equiv n\sigma_1(n)({\rm mod~}5),$
it is easily seen that $T(s)$ equals
$$\prod_{p\equiv 1({\rm mod~}5)}{1-p^{-4s}\over (1-p^{-s})(1-p^{-5s})}
\prod_{p\equiv \pm 2({\rm mod~}5)}{1-p^{-3s}\over (1-p^{-s})(1-p^{-4s})}
\prod_{p\equiv 4({\rm mod~}5)}{1\over 1-p^{-2s}}.$$
Let $\chi_c$ be the character of $(\mathbb Z/5\mathbb Z)^*$ that is determined by
$\chi_c({\bar 2})=i$ and $\chi_5$ be the character that is determined by
$\chi_5({\bar 2})=-1$.
Denote
$$\prod_{p\equiv 1({\rm mod~}5)}\left({1-p^{-4s}\over 1-p^{-5s}}\right)^4
\prod_{p\equiv \pm 2({\rm mod~}5)}{(1-p^{-3s})^4\over (1-p^{-2s})^2(1-p^{-4s})^3}
\prod_{p\equiv 4({\rm mod~}5)}(1-p^{-2s})^{-2}$$
by $H(s)$. Note that
\begin{equation}
\label{basisvergelijking}
T(s)^4=(1-5^{-s})^3H(s)\zeta(s)^3{L(s,\chi_c)L(s,{\bar \chi_c})
\over L(s,\chi_5)}.
\end{equation}
Put
$$D=\prod_{p\equiv 1({\rm mod~}5)}{1-p^{-4}\over 1-p^{-5}}
\prod_{p\equiv \pm 2({\rm mod~}5)}{1-p^{-3}\over 
(1-p^{-2})^{1/2}(1-p^{-4})^{3/4}}
\prod_{p\equiv 4({\rm mod~}5)}{1\over \sqrt{1-p^{-2}}}.$$
Thus the first order Landau-Ramanujan constant $C$ in this case, which was not
written down by Ramanujan, equals
$$C={1\over \Gamma({3\over 4})}\left({64L(1,\chi_c)L(1,{\bar \chi_c})\over 125L(1,\chi_5)}\right)^{1/4}D.$$
At \cite[p. 388]{edwards} the $L$-values above are given (for an excellent discussion of how
to compute these values see \cite[6.5]{edwards}); we have 
$$L(1,\chi_c)={\overline{L(1,{\bar \chi_c})}}
={2\pi\over 25}(3-i)\left(\sin({2\pi\over 5})+i\sin({4\pi\over 5})\right)$$ and 
$L(1,\chi_5)=\log({3+\sqrt{5}\over 2})/\sqrt{5}$.
On using that $\sin^2({2\pi/5})=(5+\sqrt{5})/8$ and $\sin^2({4\pi/5})=(5-\sqrt{5})/8$, we deduce
that
$L(1,\chi_c)L(1,{\bar \chi_c})={2\pi^2/25}$.
Alternatively we may deduce the latter equality by noting that
$${h(K)2^{r_1}(2\pi)^{r_2}R(K)\over \sqrt{|d(K)|}w(K)}={\rm Res}_{s=1}\zeta_{\mathbb Q(\zeta_5)}(s)
=L(1,\chi_c)L(1,{\bar \chi_c})L(1,\chi_5),$$
where $\zeta_{\mathbb Q(\zeta_5)}(s)$ denotes the Dedekind zeta-function of the cyclotomic
field $K:={\mathbb Q(\zeta_5)}$. It is not difficult to show that $h(K)=1$,
$r_1=0$, $r_2=2$, $R(K)=\log((3+\sqrt{5})/2)$, $d(K)=125$ and $w(K)=10$ (all of this
can be deduced from results proved e.g. in \cite{washington}).
We thus obtain that
$$C={4\over 5\Gamma({3\over 4})}\left({\pi^2\over 2\sqrt{5}\log({3+\sqrt{5}\over 2})}\right)^{1/4}D.$$
\indent Using the prime number theorem for arithmetic progressions we see that the conditions of
Lemma  \ref{edmundo} are satisfied and from (\ref{basisvergelijking}) we deduce that 
$$4B_t=-3\gamma
-2\Re\left({L'(1,{\chi_c})\over L(1,{\chi_c})}\right)
+{L'(1,\chi_5)\over L(1,\chi_5)}
-{3\over 4}\log 5 + A_{\pm 1}+A_{\pm 2}.$$
with 
$$A_{\pm 1}=\sum_{p\equiv 1({\rm mod~}5)}\log p\left({-16\over p^4-1}+{20\over p^5-1}\right)
+4\sum_{p\equiv 4({\rm mod~}5)}{\log p\over p^2-1}{\rm ~~~~~and~}$$
$$A_{\pm 2}=\sum_{p\equiv \pm 2({\rm mod~}5)}\log p\left({4\over p^2-1}-{12\over p^3-1}
+{12\over p^4-1}\right).$$

\section{Divisibility of tau by 7}
At \cite[p. 52]{bono} Ramanujan makes a statement of the form 
(\ref{valseanalogie}) with $q=7$ and $\delta=1/2$.
Put $t_n=0$ if $7|\tau(n)$ and $t_n=1$ otherwise. 
Using that $\tau(n)\equiv n\sigma_3(n)({\rm mod~}7)$ it is easily seen that
$$T(s)=\prod_{p\equiv 3,5,6({\rm mod~}7)}{1\over 1-p^{-2s}}
\prod_{p\equiv 1,2,4({\rm mod~}7)}{1-p^{-6s}\over (1-p^{-s})(1-p^{-7s})}.$$
A simple compuation shows that 
$$T(s)^2=\zeta(s)L(s,\chi_{-7})(1-7^{-s})
\prod_{p\equiv 3,5,6({\rm mod~}7)}{1\over 1-p^{-2s}}
\prod_{p\equiv 1,2,4({\rm mod~}7)}\left({1-p^{-6s}\over 1-p^{-7s}}\right)^2.$$
where $\chi_{-7}$ denotes the usual Kronecker character
of the number field $\mathbb Q(\sqrt{-7})$. From this relation
we then obtain
$$2B_t=-\gamma-{\log 7\over 6}-{L'(1,\chi_{-7})\over L(1,\chi_{-7})}
+2\sum_{p\equiv 3,5,6({\rm mod~}7)}{\log p\over p^2-1}$$
$$+\sum_{p\equiv 1,2,4({\rm mod~}7)}\log p\left({14\over p^7-1}
-{12\over p^6-1}\right).$$

\section{Divisibility of tau and lambda by 3}
\label{DRIE}
At \cite[p. 64]{bono} Ramanujan makes two statements of the form 
(\ref{valseanalogie}) with $q=3$ and $\delta=1/2$.
Put $t_n=0$ if 
$3|\tau(n)$ and
$t_n=1$ otherwise. 
Using that $\tau(n)\equiv n\sigma_1(n)({\rm mod~}3)$, where 
$\sigma_1(n)$ denotes the sum of the positive divisors of $n$, it is
easy to see that $t_n$ is multiplicative and that 
\begin{equation}
\label{doorrama}
T(s):=\sum_{n=1}^{\infty}{t(n)\over n^s}=\prod_{p\equiv 2({\rm  mod~}3)}
{1\over 1-p^{-2s}}\prod_{p\equiv 1({\rm mod~}3)}{1+p^{-s}\over 1-p^{-3s}}.
\end{equation}
A simple compuation shows that 
$$T(s)^2=\zeta(s)L(s,\chi_{-3})(1-3^{-s})
\prod_{p\equiv 2({\rm mod~}3)}{1\over 1-p^{-2s}}
\prod_{p\equiv 1({\rm mod~}3)}\left({1-p^{-2s}\over 1-p^{-3s}}\right)^2,$$
where $\chi_{-3}$ denotes the Kronecker character of the number field $\mathbb Q(\sqrt{-3})$. 
From this relation
we then obtain
$$2B_t=-\gamma-{\log 3\over 2}-{L'(1,\chi_{-3})\over L(1,\chi_{-3})}
+2\sum_{p\equiv 2({\rm mod~}3)}{\log p\over p^2-1}$$
$$+\sum_{p\equiv 1({\rm mod~}3)}\log p\left({6\over p^3-1}
-{4\over p^2-1}\right).$$
\indent Let $\lambda(n)$ denote the number of partitions of $n$ as the sum of integers which
are not multiples of $9$. 
Put $l_n=0$ if 
$3|\lambda(n)$ and
$l_n=1$ otherwise. Ramanujan shows that 
\begin{equation}
\label{koppeling}
\sum_{k\le x}l_k=\sum_{k\le 3x+1}t_k
\end{equation}
and then
states, \cite[(11.8a)]{bono}, that it can be shown by transcendental methods that
$$\sum_{k=1}^nl_k=C\int_1^n{dx\over (\log x)^{1/2}}+O\left({n\over (\log n)^r}\right)$$
and  
$$\sum_{k=1}^nt_k={C\over 3}\int_1^n{dx\over (\log x)^{1/2}}+O\left({n\over (\log n)^r}\right),
$$
where he gives an explicit expression for $C$. Though these claims are correct for $r\le 3/2$, 
{\it a priori} they cannot be both true for $r>3/2$, as we then trivially have 
from (\ref{koppeling}) that $C_2(l)=C_2(t)-{1\over 2}\log 3$, 
whereas the truth of both claims of Ramanujan would imply that $C_2(l)=C_2(t)$. 

\section{Divisibility of tau by 691}
At \cite[p. 66]{bono} Ramanujan makes a statement of the form 
(\ref{valseanalogie}) with $q=691$ and $\delta=1/690$. 
He did not write down an explicit first order constant.
The truth
of this assertion for $r\le 691/690$ was 
first established by G.N. Watson \cite{watson1}. In this note
it will be shown, however,
that the statement is false for every $r>691/690$.\\
\indent It is not difficult to show, as Ramanujan did, that 
$\tau(n)\equiv \sigma_{11}(n) ({\rm mod~}691)$.  
Let $\nu(p)$ be the smallest integer $>1$, such that $p^{\nu(p)}\equiv 1({\rm mod~}691)$.
We put $\nu(691)=\infty$. Note that $\sigma_{11}(p^k)\equiv 0({\rm mod~}691)$ if and 
only if $k\equiv \nu(p)-1({\rm mod~}\nu(p))$. In case $p=691$ we interpret this congruence
as never being satisfied and $1-691^{-\nu(691)s}$ as being $1$. We thus can write
$$T(s)=\prod_{p}{1-p^{-(\nu(p)-1)s}\over (1-p^{-s})(1-p^{-\nu(p)s})}.$$
Around $s=1$ this function is quite close to $\zeta(s)$ 
(we have 
$T(s)=\sum_{n=1}^{11053}n^{-s}-\sum_{m=1}^{8}(m1381)^{-s}-5527^{-s}-8291^{-s}+\sum_{n=11054}^{\infty}t_n/n^s$) and hence we expect $B_t$ to be
close to $-\gamma$, cf. Example 1, which is indeed the case by Table 1.\\
\indent Notice that each local factor of $T(s)$ has the term $1-p^{-s}$ in it, unless
$p\equiv -1({\rm mod~}691)$, in which case the local factor is $(1-p^{-2s})^{-1}$. By multiplying
$T(s)$  with $\prod_{p\equiv -1({\rm mod~}691)}(1-p^{-s})^{-1}$ we can then write it as
$\zeta(s)h(s)$ with $h(s)$ a regular function for Re$(s)>1/2$. On noting that the
primes $p$ with $p\equiv -1({\rm mod~}691)$  are exactly the primes that split completely
in $\mathbb Q(\zeta_{691})$ but not in $\mathbb Q(\cos({2\pi\over 691}))$, the product
$\prod_{p\equiv -1({\rm mod~}691)}(1-p^{-s})^{-1}$ can be expressed in terms of the
Dedekind zetafunctions of the latter two fields and some regular function for Re$(s)>1/2$,
which can be explicitly determined using the splitting behaviour of a prime $p$ in these
fields. Using the factorisation of the these Dedekind zetafunctions in terms of L-series,
we then obtain the following identity, with $\chi_c$ the character uniquely determined
by $\chi_c({\bar 3})=\exp(2\pi i/690)$;
\begin{eqnarray}
T(s)^{690}&=&\zeta(s)^{689}(1-691^{-s})^{-1}
L(s,\chi_c)\prod_{j=1}^{344}{L(s,\chi^{2j+1}_c)
\over L(s,\chi^{2j}_c)}\nonumber\cr
&&\prod_{p\equiv -1({\rm mod~}691)}(1-p^{-2s})^{-345}
\prod_{p\equiv 1({\rm mod~691})}\left({1-p^{-690s}\over 1-p^{-691s}}\right)^{690}\nonumber\cr
&&\prod_{p\ne 691\atop 2|\nu(p),~\nu(p)\ge 4}
\left({1+p^{-\nu(p){s\over 2}}\over 1-p^{-\nu(p){s\over 2}}}\right)^{690\over \nu(p)}
\prod_{p\atop 2<\nu(p)<691}\left({1-p^{-(\nu(p)-1)s}\over 1-p^{-\nu(p)s}}\right)^{690},
\nonumber
\end{eqnarray}
the truth of which is most easily established by checking that the local factors on both
sides agree for every prime $p$. As before a formula for $B_t$ can now be easily written
down, but
for reasons of space we leave this to the interested reader. It turns out that the
contribution of the last four products in the formula for $T(s)^{690}$ to $B_t$ is less
than $10^{-5}$ in absolute value. We thus have
\begin{equation}
\label{zuchtsteun}
B_t\approx {\log 691\over 690^2}-{689\over 690}\gamma-{1\over 690}\sum_{j=0}^{344}
{L'(1,\chi_c^{2j+1})\over L(1,\chi_c^{2j+1})}+{1\over 690}\sum_{j=1}^{344}
{L'(1,\chi_c^{2j})\over L(1,\chi_c^{2j})}.
\end{equation}
with an error of at most $10^{-5}$.

\section{Divisibility of tau by 23}
At \cite[p. 80]{bono} Ramanujan makes a statement of the form 
(\ref{valseanalogie}) with $q=23$ and $\delta=1/2$. Using a trick of Wilton 
\cite{wilton}, $T(s)$ can be
found in this case more easily than by Ramanujan's approach. Using Euler's identity we note
that
\begin{eqnarray}
\sum_{n=1}^{\infty}\tau(n)x^n&=&x\{(1-x)(1-x^2)\cdots\}^{23}(1-x)(1-x^2)\cdots\nonumber\\
&\equiv&xk(x^{23})\sum_{m=-\infty}^{\infty}(-1)^mx^{m(3m+1)\over 2}
~({\rm mod~}23),\nonumber
\end{eqnarray}
where $k(x)=\prod_{n=1}^{\infty}(1-x^n)$.
Now if we also apply Euler's identity to
$k(x^{23})$, then we obtain
$$\sum_{n=1}^{\infty}\tau(n)x^n\equiv x\sum_{r=-\infty}^{\infty}\sum_{m=-\infty}^{\infty}
(-1)^{m+r}x^{{m(3m+1)\over 2}+23{r(r+1)\over 2}}~({\rm mod~}23).$$
From the latter identity Wilton's congruences are easily deduced;
$$
\tau(p)\equiv
\cases{1~({\rm mod~}23) &if $p=23$;\cr
0~({\rm mod~}23) &if $(p/23)=-1$;\cr
2~({\rm mod~}23) &if $p=U^2+23V^2$ with $U\ne 0$;\cr
-1~({\rm mod~}23) &for other $p\ne 23$.}
$$
Using that $\tau(p^{k+1})=\tau(p)\tau(p^k)-p^{11}\tau(p^{k-1})$ for $k\ge 1$, we can now easily
compute $\tau(p^k)$ modulo $23$.
Let ${\cal S}_1$ denote the set of primes $p$ such that $p$ is a quadratic non-residue mod
$23$. Let ${\cal S}_3$ denote the set of primes $p$ which can be written as $U^2+23V^2$ with $U\ne 0$.
Let ${\cal S}_2$ be the set of remaining primes $p\ne 23$. 
(The sets ${\cal S}_1,{\cal S}_2,{\cal S}_3$ have natural densities of 
respectively, $1/2$, $1/3$ and $1/6
$, which can be shown
using that the class number of $\mathbb Q(\sqrt{-23})$ equals 3.)
We now find that
$$T(s)={1\over 1-23^{-s}}\prod_{p\in {\cal S}_1}{1\over 1-p^{-2s}}
\prod_{p\in {\cal S}_2}{1+p^{-s}\over 1-p^{-3s}}
\prod_{p\in {\cal S}_3}{1-p^{-22s}\over (1-p^{-s})(1-p^{-23s})}.$$
The factor at 23 is not present in Ramanujan's formula (17.6), although it should be
there according to his argument. This leads then to an incorrect formula for the 
first order constant (at the bottom of \cite[p. 80]{bono}). On should replace the factor
$66^{1/2}23^{-3/4}$ by $23^{1/4}\sqrt{3/22}$. It is also clear that where
Ramanujan writes `all primes of the form $23a^2+b^2$, he excludes the prime 23.\\
\indent We have
$$T(s)^2={\zeta(s)L(s,\chi_{-23})\over 1-23^{-s}}\prod_{p\in {\cal S}_1}
{1\over 1-p^{-2s}}\prod_{p\in {\cal S}_2}
\left({1-p^{-2s}\over 1-p^{-3s}}\right)^2
\prod_{p\in {\cal S}_3}
\left({1-p^{-22s}\over 1-p^{-23s}}\right)^2.
$$
From this we easily deduce that
\begin{eqnarray}
B_t&=&-{\gamma\over 2}-{L'(1,\chi_{-23})\over 2L(1,\chi_{-23})}+{\log 23\over 44}\nonumber\cr
&&+\sum_{p^{11}\equiv -1({\rm mod~}23)}{\log p\over p^2-1}+
\sum_{p^{11}\equiv 1({\rm mod~}23)}\log p\left({3\over p^3-1}-{2\over p^2-1}\right)\nonumber\cr
&&
+\sum_{p=U^2+23V^2\atop p>23}\log p\left({2\over p^2-1}-{3\over p^3-1}+{23\over p^{23}-1}
-{22\over p^{22}-1}\right).\nonumber
\end{eqnarray}
Remark. To the reader familiar with Cox's beautiful book \cite{cox}, we suggest as an exercise
showing that $p\in {\cal S}_3$ if and only if $\left({p\atop 23}\right)=1$ and
the congruence $x^3\equiv x+1({\rm mod~}p)$ has an integer solution.
\section{Numerical evaluation of the second order constants}
\label{numerical}
The expressions obtained for the various $B_t$ involve both prime sums and
values of $L$ and $L'$ at $s=1$.
The prime sums we evaluate termwise and estimate the tail using that, for $k>1$
and $x\ge 7481$,
$$\sum_{p>x}{\log p\over p^k-1}\le {x\over x^k-1}(-0.98+1.017{k\over k-1}),$$
which follows easily on using the estimate $0.98x\le \theta(x)\le 1.017x$ for
$x\ge 7481$ \cite{RS}.\\
\indent The $L$ and $L'$ values above can be evaluated using generalized Euler
constants for arithmetical progressions. We define
$$\gamma_k(r,m)
:=\lim_{x\rightarrow \infty}
\left\{\sum_{0<n\le x\atop n\equiv r({\rm mod~}m)}{\log ^k n\over n}-
{\log^{k+1}x\over m(k+1)}\right\}.$$
Note that $\gamma_0(0,1)=\gamma$, Euler's constant.
Let $\chi$ be a non-principal character modulo $m$. It is not difficult 
to show \cite{knopfmacher} that for $k\ge 0$ we have
$$L^{(k)}(1,\chi)=(-1)^k\sum_{r=1}^m \chi(r)\gamma_k(r,m).$$
Using Proposition 12 of \cite{dilcher}, the Euler constants
$\gamma_k(r,m)$ can be computed with any degree of precision and
thus the same holds true for $L^{(k)}(1,\chi)$.\\
\indent For $q=5$ we find, using Dilcher's Table 1 \cite[S21]{dilcher}, that
$L'(1,\chi_5)/L(1,\chi_5)=0.82767947\cdots$ and
$L'(1,\chi_c)/L(1,\chi_c)=0.15786453\cdots-i0.08833613\cdots$.
For $q=7$ we find using Dirichlet's formula that $L(1,\chi_{-7})=\pi/\sqrt{7}$ and,
 using Dilcher's Table 1, that $L'(1,\chi_{-7})=0.01856598\cdots$. 
 The quotient $L'(1,\chi_{-3})/L(1,\chi_{-3})$ is evaluated with many decimal accuracy
in, e.g., \cite{moree}. We have $L(1,\chi_{-23})=3\pi/\sqrt{23}$ and, on implementing
Proposition 12 of \cite{dilcher} in Maple, we find $L'(1,\chi_{-23})=-0.82955295\cdots$.
Similarly we find that the sum involving the odd, respectively
even characters in (\ref{zuchtsteun}) equal $1.9018228\cdots$, respectively 
$5.10942407\cdots$ (note that {\it a priori}
these sums must be real).\\
\indent For $q=3$ we can use the relationship $-\zeta'(2)/\zeta(2)=\sum_{p}{\log p/(p^2-1)}$ to
rewrite $B_t$ as
$$2B_t=6\sum_{p\equiv 2({\rm mod~}3)}{\log p\over p^2-1}+4{\zeta'(2)\over \zeta(2)}
-{L'(1,\chi_{-3})\over L(1,\chi_{-3})}-\gamma+6\sum_{p\equiv 1({\rm mod~}3)}{\log p\over p^3-1}.$$
Since all but the last term were either computed with high accuracy in \cite{moree} or are
easily computable with high accuracy in Maple, we now obtain that
$B_t=-0.5349219\cdots$.\\
\indent Let us for a function $f$ satisfying the conditions of Theorem \ref{oud}
define $H_f(x):=\sum_{n\le x}\Lambda_f(n)/n-\tau \log x$. 
The function $\Lambda_f$ is most easily computed by computing minus the logarithm of
the generating series of $f$.
The numbers $H_f(10^5)$
and $H_f(10^6)$ ought to be approximations of $B_f$. 
The function $b$ is the indicator function of the set of integers that can be written
as a sum of two squares. This is, as was already known to Fermat, a multiplicative function.
The final column in Table 1
gives Ramanujan's predicted value for $C_2$.\\

\centerline{{\bf Table 1:} Numerical values}
\begin{center}
\begin{tabular}{|c|c|c|c|c|c|c|}\hline
&$f$&$H_f(10^5)$&$H_f(10^6)$&$B_f$&$C_2(f)$&R's $C_2$\\ \hline\hline
$B(x)$&b&$+0.163\cdots$&$+0.162\cdots$&$+0.1638\cdots$&$0.5819\cdots$&1/2\\ \hline
$q=5$&t&$-0.401\cdots$&$-0.400\cdots$&$-0.3995\cdots$&$0.1501\cdots$&1/4\\ \hline
$q=7$&t&$-0.232\cdots$&$-0.232\cdots$&$-0.2316\cdots$&$0.3841\cdots$&1/2\\ \hline
$q=3$&t&$-0.532\cdots$&$-0.534\cdots$&$-0.5349\cdots$&$0.2325\cdots$&1/2\\ \hline
$q=691$&t&$-0.571\cdots$&$-0.571\cdots$&$-0.5717\cdots$&$0.0006\cdots$&1/690\\ \hline
$q=23$&t&$-0.217\cdots$&$-0.217\cdots$&$-0.2166\cdots$&$0.6083\cdots$&1/2\\ \hline
\end{tabular}
\end{center}

\noindent Remark. The computations were not carried out far enough to determine the fifth
digit in $-0.2166\cdots$; it is either a $6$ or a $7$.\\

\noindent Now we are in the position to prove the following result.
\begin{Thm}
All assertions made by Ramanujan in his `unpublished' manuscript on the 
partition and
tau functions {\rm \cite{bono}} of the format {\rm (\ref{valseanalogie})} are false.
\end{Thm}
{\it Proof}. Let $q$ be a prime from Table 1. Assume (\ref{valseanalogie}) holds
true for $r>1+\delta$, with $\delta$ as in the last column of Table 1. Then the
second order coefficient equals $\delta$, which does not match the value of $C_2(t)$
given in Table 1. For the function $\lambda$ from Section \ref{DRIE} we have
$C_2(l)=0.2325\cdots-{1\over 2}\log 3\ne 0.5$. \qed

\section{On a 1928 paper of Geraldine Stanley}
\label{Geraldine}

The purpose of Stanley's paper \cite{stanley} is to show that two assertions due to Ramanujan
are false. The first assertion was already mentioned in the introduction.
Stanley's analysis of this case contains, unfortunately, several misprints
and errors, which are corrected in \cite{shanks}.\\
\indent The second assertion concerns the $5|\tau (n)$ case.
In one of the footnotes we read:
"In discussing this question I have used a manuscript of Prof. Hardy, who
at one time intended to complete Ramanujan's work".
Hardy made some headway with this and then apparently later asked Stanley
to fill in the further details.
The purpose of this section is to correct Geraldine Stanley's analysis of this case and
point out typo's. With respect to the analysis of the first assertion Shanks 
\cite[p. 75]{shanks} has written: "In fact, there
are several errors, and these nullify the 
proof that Ramanujan's second term is wrong".
Indeed, we will see that the same wording applies to her analysis of the $5|\tau (n)$ case.
For information on Geraldine Stanley, the reader is referred to
\cite{rankin}.\\
\indent  In the definition of $\psi_1(s)$, for $p_2^{-s}$ read $p_1^{-s}$.
For $p_2^{-s}$ read $p_1^{-5s}$. A few lines down read $(1,i,-i,-1)$ instead
of $(1,i,-i,1)$. 
For $(1-5^{-s})^{5/4}$ in the formula for $h(s)$ read
$(1-5^{-s})^{3/4}$. 
In the denominator of $h(s)$, the factor
$(1-p_3^{-4s})^{3/4}$ has to be added in the denominator (the so corrected $h(s)$ is thus invariant under
permuting $p_2$ and $p_3$). 
In the formula for $A$ at the bottom of p. 236, the factor
$(4/5)^{3/2}$ has to be changed to $4/5$. This formula shows that the omission of
the factor $(1-p_3^{-4s})^{3/4}$ in the formula for $h(s)$ was not a mathematical mistake.
Also the formula for $a_1/A$ given at p. 237 shows that this factor was in the original formula
for $h(s)$. The exponent $5/4$ instead of $3/4$ is consistenly worked with in the remainder
of the paper though. This leads then to $(4/5)^{3/2}$ instead of $4/5$ in the formula for $A$ and
to $5(\log 5)/16$ instead of $3(\log 5)/16$ in the formula for $a_1/A$.\\
\indent From the numerical point of view the formula for $A$ is quite awkward since the
values $L_j(1)$, $2\le j\le 4$ have been expanded as Euler products and have been regrouped.
(On not doing this a numerically more convenient expression for $A$ can be obtained, cf.
Section \ref{five}.)
Thus we have to read 
$$\prod\left({1-p_4^{-1}\over (1-p^{-1})^{1/4}}\right){\rm ~as~}
\lim_{x\rightarrow \infty}{\prod_{p_4\le x}(1-p_4^{-1})\over 
\prod_{p\le x}(1-p^{-1})^{1/4}}.$$
Indeed, by the Mertens' theorem for arithmetic progressions \cite{williams} this limit exists. 
\indent The factor $1/(s+1)^{3/4}$ at the middle of p. 237 has to be replaced by $1/(s-1)^{3/4}$.
The formula for $b_1$ at p. 237 is off by a minus sign (cf. \cite[pp. 720-723]{watson2}).
It should read $b_1=-a_1\Gamma({5\over 4})/\pi\sqrt{2}$.
Consequently the term $(a_1-A)$ in the first formula for $T(x)$ becomes $(-a_1-A)$ and
in the formula for $a_1/A$, $a_1/A$ has to be replaced by $-a_1/A$. This sign error nullifies the
last two sentences of the paper. (Note that apart from the sign error,
$5(\log 5)/16$ has to be replaced by $3(\log 5)/16$ and
$L_4(1)$ by $4L_4(1)$.)

\end{document}